\documentclass[12pt, a4paper]{amsart}
\usepackage[utf8]{inputenc}
\usepackage{amsxtra,amssymb,amsthm,amsmath,amscd,mathrsfs, epsfig, eufrak}
\usepackage{amscd, amsmath, mathrsfs, amssymb, amsthm, amsxtra, bbding, epsfig, eucal, eufrak, graphicx, latexsym, mathrsfs, mathbbol, bbold}
\usepackage[all]{xy}

\usepackage{tikz}

\usepackage[normalem]{ulem}
\usepackage{soul}
\usepackage{color}

\setstcolor{red}

\makeatletter
\@namedef{subjclassname@2010}{%
\textup{2010} Mathematics Subject Classification}
\makeatother

\theoremstyle{plain}
\newtheorem{theorem}{Theorem}
\newtheorem{proposition}{Proposition}[section]
\newtheorem{lemma}[proposition]{Lemma}
\newtheorem{corollary}{Corollary}

\theoremstyle{remark}

\numberwithin{equation}{section}

\addtocounter{footnote}{1}

\begin{document}


\title[\tiny On a kind of generilized multi-harmonic sums]
{On a kind of generalized multiharmonic sums}
\author[\tiny Jiaqi Wang \& Rong Ma]{Jiaqi Wang \& Rong Ma}

\address{%
Jiaqi Wang
\\
School of Mathematics and Statistics
\\
Northwestern Polytechnical University
\\
Xi'an
\\
Shaanxi 710072
\\
China}
\email{1727561320@qq.com}

\address{%
Rong Ma
\\
School of Mathematics and Statistics
\\
Northwestern Polytechnical University
\\
Xi'an
\\
Shaanxi 710072
\\
China}
\email{marong@nwpu.edu.cn}

\date{\today}

\begin{abstract}
Let $p$ be an odd prime, Jianqiang Zhao has established a curious congruence, which is
$$
 \sum_{i+j+k=p \atop i,j,k > 0} \frac{1}{ijk} \equiv -2B_{p-3}\pmod p ,
$$
where $B_{n}$ denotes the $n$-th Bernoulli number. In this paper, we will generalize this kind of sums and prove a family of similar congruences modulo prime powers $p^r$.
\end{abstract}

\footnote {2020 Mathematics Subject Classification{: Primary 11B68; Secondary 11A07.}}

\maketitle

\section{Introduction}
For any prime $p \geq 3$ , $B_n$ is the $n$-th Bernoulli number. In 2007, Jianqiang Zhao \textsuperscript{\cite{JQZ2007}} first established a curious congruence
\begin{equation}\label{PNT:weak}
    \sum_{i+j+k=p \atop i,j,k > 0} \frac{1}{ijk} \equiv -2B_{p-3} \pmod p.
\end{equation}

In recent years, more and more researchers have generalized the formula above. For example, Xia Zhou and Tianxin Cai \textsuperscript{\cite{XZ2007}} used mathematical induction and the difference method and generalized the formula
(1.1). For any prime $p \geq 5$ and $n \leq p-2$, they got
\begin{equation}\label{PNT:strong}
\sum_{l_{1}+l_{2}+...+l_{n}=p \atop l_{1},l_{2},...,l_{n}>0} \frac{1}{l_{1}l_{2}...l_{n}} \equiv \left\{
    \begin{aligned}
        -(n-1)!B_{p-n} \pmod p ,\quad  \mbox{if} \quad 2 \nmid n;\\
         -\frac{n}{2(n+1)}n!B_{p-n-1}p \pmod {p^2} ,\quad \mbox{if} \quad 2 \mid n .\\
    \end{aligned}
\right
.
\end{equation}

In addition, researchers have generalized formula (1.1) in various aspects. For example, in 2014, Liuquan Wang and Tianxin Cai \textsuperscript{\cite{Lq2014}} got a generalized curious congruence. For any prime $p \geq 3$ and any positive integer $r$, they obtained
\begin{equation}\label{PNT:RH}
    \sum_{i+j+k=p^r \atop i,j,k \in P_{p}} \frac{1}{ijk} \equiv -2p^{r-1}B_{p-3} \pmod{p^r},
\end{equation}
where $P_p$ denotes the positive integers which are coprime to $p$.

In 2015, Liuquan Wang \textsuperscript{\cite{Lq2015}} used the inclusion-exclusion principle and congruences of multiple harmonic sums and generalized the result, that is
\begin{equation}\label{def:Sfx}
    \sum_{i_1+i_2+...+i_5=p^r \atop i_1,i_2,...,i_5 \in P_{p}} \frac{1}{i_1i_2...i_5} \equiv -\frac{5!}{6}p^{r-1}B_{p-5}\pmod{p^r}.
\end{equation}

Subsequently, in 2021, Liuquan Wang and Jianqiang Zhao\textsuperscript{\cite{LqJQ}} have generalized once again and got the following congruence
\begin{equation}
    \sum_{l_1+...+l_7=p^r \atop l_1,...,l_7 \in P_p}\frac{1}{l_1...l_7}\equiv-\frac{7!}{10}p^{r-1}B_{p-7} \pmod{p^r}.
\end{equation}
For any positive integer $n$, we define $Z(n)$ as
\begin{equation}
    Z(n):=\sum_{i+j+k=n \atop i,j,k \in P_n}\frac{1}{ijk},
\end{equation}
where $P_n$ denotes the positive integers which are coprime to $n$.

Although there have been numerous generalizations about formula (1.1), it remains many problems which have not yet to be conducted. In this paper, we shall further generalize the congruence (1.3) from $i+j+k=p^r$ to $i+j+k=p_1^{r_1}...p_n^{r_n}$ where $p_1,...,p_n$ are distinct odd primes and $r_1,...,r_n$ are positive integers. We will also get the following theorems.

\begin{theorem}\label{thm1}
Let $p_1,...,p_n$ be distinct odd primes, then we have
\begin{equation}
\begin{aligned}
    Z(p_1...p_n)&\equiv-2p_2...p_nB_{p_1-3}\left(1+\sum_{l_1=2}^{n}\frac{2}{p_{l_1}^4}-\sum_{l_1,l_2 \in \{2,...,n\} \atop l_1 \neq l_2}\frac{2}{p_{l_1}^4p_{l_2}^4}+...+(-1)^n\frac{2}{p_{2}^4...p_{n}^4}  \right)\\
    &+6p_2...p_nB_{p_1-3}\left(\sum_{l_1=2}^{n}\frac{2p_{l_1}^2+1}{3p_{l_1}^3}-\sum_{l_1,l_2 \in \{2,...,n\} \atop l_1 \neq l_2}\frac{2p_{l_1}^2p_{l_2}^2+1}{3p_{l_1}^3p_{l_2}^3}+...\right.\\
    &\left.+(-1)^n\frac{2p_{2}^2...p_{n}^2+1}{3p_{2}^3...p_{n}^3} \right)\pmod{p_1}.
\end{aligned}
\end{equation}
\end{theorem}

\begin{theorem}\label{thm1}
Let $p_1,...,p_n$ be distinct odd primes and $r_1,...,r_n$ be positive integers, then we have
\begin{equation}
\begin{aligned}
    Z(p_1^{r_1}...p_n^{r_n})&\equiv-2p_1^{r_1-1}p_2^{r_2}...p_n^{r_n}B_{p_1-3}\left(1+\sum_{l_1=2}^{n}\frac{2}{p_{l_1}^4}-\sum_{l_1,l_2 \in \{2,...,n\} \atop l_1 \neq l_2}\frac{2}{p_{l_1}^4p_{l_2}^4}+...+(-1)^n\frac{2}{p_{2}^4...p_{n}^4}  \right)\\
    &+6p_1^{r_1-1}p_2^{r_2}...p_n^{r_n}B_{p_1-3}\left(\sum_{l_1=2}^{n}\frac{2p_{l_1}^2+1}{3p_{l_1}^3}-\sum_{l_1,l_2 \in \{2,...,n\} \atop l_1 \neq l_2}\frac{2p_{l_1}^2p_{l_2}^2+1}{3p_{l_1}^3p_{l_2}^3}+...\right.\\
    &\left.+(-1)^n\frac{2p_{2}^2...p_{n}^2+1}{3p_{2}^3...p_{n}^3} \right)\pmod{p_1^{r_1}}.
\end{aligned}
\end{equation}
\end{theorem}
In Theorem 2, let $n=1$ and $n=2$, we can obtain previous results, see Ref.\cite{Lq2014} and \cite{Cai2016}, that is the following corollaries.
\begin{corollary}\textsuperscript{\cite{Lq2014}}
    Let $p$ be distinct odd primes and $r$ be positive integers, then we have
    \begin{equation}
    \begin{aligned}
        Z(p^{r})\equiv&-2p^{r-1}B_{p-3}\pmod{p^{r}}.
    \end{aligned}
    \end{equation}
\end{corollary}
\begin{corollary}\textsuperscript{\cite{Cai2016}}
    Let $p_1,p_2$ be distinct odd primes and $r_1,r_2$ be positive integers, then we have
    \begin{equation}
    \begin{aligned}
        Z(p_1^{r_1}p_2^{r_2})\equiv&-2p_1^{r_1-1}p_2^{r_2}B_{p_1-3}\left(1+\frac{2}{p_2^4}\right)+6p_1^{r_1-1}p_2^{r_2}\frac{2p_2^2+1}{3p_2^3}B_{p_1-3}\\
        \equiv&2(2-p_2)\left(1-\frac{1}{p_2^3}\right)p_1^{r_1-1}p_2^{r_2-1}B_{p_1-3}\pmod{p_1^{r_1}}.
    \end{aligned}
    \end{equation}
\end{corollary}

For convenience, given positive integers $n$, $m$, and any odd prime $p$, we define
$$
S(m,n;p):=\sum_{a=1 \atop (a,p)=1}^{n-1}\frac{1}{a^2}\sum_{i=1 \atop (i,p)=1}^{am-1}\frac{1}{i}
$$
and
$$
T(m,n;p):=\sum_{a=1 \atop (a,p)=1}^{n-1}\frac{1}{a^2}\sum_{i=1 \atop i \equiv sm \pmod{p}}^{am-1}\frac{1}{i}.
$$

\vskip 8mm

\section{Lemmas}
In order to prove our theorems, we need the following lemmas.
\begin{lemma}
    Let p be an odd prime and l, m, s be positive integers, then we have
    \begin{equation*}
            \sum_{k=1 \atop (k,p)=1}^{mp^l-1} \frac{1}{k^s}\equiv
\begin{cases}
 \begin{aligned}
     &0 \pmod{p^{2l-1}}, &&\text{{\rm for odd $s$ with $p-1 \mid s+1$ and $p \nmid s$},}  \\
     &0  \pmod{p^{2l}}, &&\text{{\rm for odd $s$ with $p-1 \nmid s+1$ or $p \mid s$},} \\
     &0 \pmod{p^{l-1}},&&\text{{\rm for even s with $p-1 \mid s$},}   \\
     &0 \pmod{p^{l}},&&\text{{\rm for even s with $p-1 \nmid s$}.}
 \end{aligned}
\end{cases}
    \end{equation*}
    Proof: By \cite{Cai2016}, we have the following formulas
    \begin{equation}
     \sum_{k=1 \atop (k,p)=1}^{p^l-1} \frac{1}{k^s}\equiv
\begin{cases}
 \begin{aligned}
     &0 \pmod{p^{2l-1}}, &&\text{{\rm for odd $s$ with $p-1 \mid s+1$ and $p \nmid s$},}  \\
     &0  \pmod{p^{2l}}, &&\text{{\rm for odd $s$ with $p-1 \nmid s+1$ or $p \mid s$},} \\
     &0 \pmod{p^{l-1}},&&\text{{\rm for even s with $p-1 \mid s$},}   \\
     &0 \pmod{p^{l}},&&\text{{\rm for even s with $p-1 \nmid s$}.}
 \end{aligned}
\end{cases}
    \end{equation}
    Then for any positive integer n, we conclude the formula
    \begin{equation*}
        \sum_{k=np^l+1 \atop (k,p)=1}^{(n+1)p^l-1} \frac{1}{k^s}=\sum_{k=1 \atop (k,p)=1}^{p^l-1} \frac{1}{(k+np^l)^s}\equiv \sum_{k=1 \atop (k,p)=1}^{p^l-1}\frac{1}{k^s}-snp^l\sum_{k=1 \atop (k,p)=1}^{p^l-1}\frac{1}{k^{s+1}}\pmod{p^{2l}}.
    \end{equation*}
    When s is an odd integer which satisfies $p-1\mid s+1$ and $p \nmid s$, from (2.1), we have
    \begin{equation*}
        \sum_{k=1 \atop (k,p)=1}^{p^l-1}\frac{1}{k^s}\equiv0\pmod{p^{2l-1}},\sum_{k=1 \atop (k,p)=1}^{p^l-1}\frac{1}{k^{s+1}}\equiv0\pmod{p^{l-1}}.
    \end{equation*}
    So we have,
    \begin{equation}
        \sum_{k=np^l+1 \atop (k,p)=1}^{(n+1)p^l-1} \frac{1}{k^s} \equiv0 \pmod{p^{2l-1}}.
    \end{equation}
    Similarly, when s is an odd integer which satisfies $p-1\nmid s+1$ or $p \mid s$, by (2.1), we have
    \begin{equation}
         \sum_{k=np^l+1 \atop (k,p)=1}^{(n+1)p^l-1} \frac{1}{k^s} \equiv0 \pmod{p^{2l}}.
    \end{equation}
    When s is an even integer which satisfies $p-1\mid s$, we have,
    \begin{equation}
         \sum_{k=np^l+1 \atop (k,p)=1}^{(n+1)p^l-1} \frac{1}{k^s} \equiv0 \pmod{p^{l-1}}.
    \end{equation}
    When s is an even integer which satisfies $p-1\nmid s$, we have,
    \begin{equation}
        \sum_{k=np^l+1 \atop (k,p)=1}^{(n+1)p^l-1} \frac{1}{k^s} \equiv0 \pmod{p^{l}}.
    \end{equation}
    Therefore, we get
    \begin{equation}
     \sum_{k=np^l+1 \atop (k,p)=1}^{(n+1)p^l-1} \frac{1}{k^s}\equiv
\begin{cases}
 \begin{aligned}
     &0 \pmod{p^{2l-1}}, &&\text{{\rm for odd $s$ with $p-1 \mid s+1$ and $p \nmid s$},}  \\
     &0  \pmod{p^{2l}}, &&\text{{\rm for odd $s$ with $p-1 \nmid s+1$ or $p \mid s$},} \\
     &0 \pmod{p^{l-1}},&&\text{{\rm for even s with $p-1 \mid s$},}   \\
     &0 \pmod{p^{l}},&&\text{{\rm for even s with $p-1 \nmid s$}.}
 \end{aligned}
\end{cases}
    \end{equation}
    By the method of summation, we have finally got
      \begin{equation*}
            \sum_{k=1 \atop (k,p)=1}^{mp^l-1} \frac{1}{k^s}\equiv
\begin{cases}
 \begin{aligned}
     &0 \pmod{p^{2l-1}}, &&\text{{\rm for odd $s$ with $p-1 \mid s+1$ and $p \nmid s$},}  \\
     &0  \pmod{p^{2l}}, &&\text{{\rm for odd $s$ with $p-1 \nmid s+1$ or $p \mid s$},} \\
     &0 \pmod{p^{l-1}},&&\text{{\rm for even s with $p-1 \mid s$},}   \\
     &0 \pmod{p^{l}},&&\text{{\rm for even s with $p-1 \nmid s$}.}
 \end{aligned}
\end{cases}
    \end{equation*}
\end{lemma}

\begin{lemma}
    Let p be an odd prime, n and m be positive integers coprime to p, then we have
    \begin{equation*}
        S(m,np;p)\equiv nm^2B_{p-3}\pmod{p}.
    \end{equation*}
    Proof: For any positive integer n and r, we have
    \begin{equation}
        \sum_{a=1}^{n-1}a^r=\frac{1}{r+1}\sum_{k=0}^{r}\binom{r+1}{k}B_kn^{r+1-k},
    \end{equation}
    hence
    \begin{equation*}
    \begin{aligned}
         \sum_{a=1 \atop (a,p)=1}^{np-1}\frac{1}{a^2}\sum_{i=1 \atop (i,p)=1}^{am-1}\frac{1}{i}&\equiv\sum_{a=1 \atop (a,p)=1}^{np-1}\frac{1}{a^2}\sum_{i=1 \atop (i,p)=1}^{am-1}i^{p-2}\equiv\sum_{a=1 \atop (a,p)=1}^{np-1}\frac{1}{a^2}\sum_{i=1}^{am-1}i^{p-2}\\
         &\equiv\sum_{a=1 \atop (a,p)=1}^{np-1}\frac{1}{a^2}\frac{1}{p-1}\sum_{k=0}^{p-2}\binom{p-1}{k}B_k(am)^{p-1-k}\\
         &\equiv\frac{m^2}{p-1}\sum_{k=0}^{p-2}\binom{p-1}{k}B_k\sum_{a=1 \atop (a,p)=1}^{np-1}(am)^{p-3-k}\pmod{p}.\\
    \end{aligned}
    \end{equation*}
    When $k \neq p-3$, we have
    \begin{equation*}
         \sum_{a=1 \atop (a,p)=1}^{np-1}(am)^{p-3-k} \equiv0\pmod{p}.
    \end{equation*}
    So we have
    \begin{equation*}
       S(n,np;p)\equiv\frac{m^2}{p-1}\binom{p-1}{p-3}B_{p-3}n(p-1)\equiv nm^2B_{p-3} \pmod{p}.
    \end{equation*}
\end{lemma}

\begin{lemma}
    Let p be an odd prime, n, m be positive integers coprime to p and $\alpha \geq 2$ positive integer, then we have
    \begin{equation}
        S(m,np^{\alpha};p)\equiv np^{\alpha-1}m^2B_{p-3}\pmod{p^{\alpha}}
    \end{equation}
    Proof: Let $a=s+np^{\alpha-1}t$, $1\leq s \leq np^{\alpha-1}-1$, $1 \leq t \leq p-1$, $(s,p)=1$, then
    \begin{equation*}
        \begin{aligned}
            S(np^{\alpha};p)&=\sum_{t=0}^{p-1}\sum_{s=1 \atop (s,p)=1}^{np^{\alpha-1}-1}\frac{1}{(s+np^{\alpha-1}t)^2}\sum_{i=1 \atop (i,p)=1}^{(s+np^{\alpha-1}t)m-1}\frac{1}{i}\\
            &\equiv \sum_{t=0}^{p-1}\sum_{s=1 \atop (s,p)=1}^{np^{\alpha-1}-1}\frac{1}{s^2}(1-\frac{2np^{\alpha-1}t}{s})\sum_{i=1 \atop (i,p)=1}^{sm-1}\frac{1}{i}\\
            &+\sum_{t=0}^{p-1}\sum_{s=1 \atop (s,p)=1}^{np^{\alpha-1}-1}\frac{1}{s^2}(1-\frac{2np^{\alpha-1}t}{s})\sum_{i=sm \atop (i,p)=1}^{(s+np^{\alpha-1}t)m-1}\frac{1}{i}.\\
        \end{aligned}
    \end{equation*}
    Trivially, we have
    \begin{equation*}
        2np^{\alpha-1}\sum_{t=0}^{p-1}t\sum_{s=1 \atop (s,p)=1}^{np^{\alpha-1}-1}\frac{1}{s^3}\sum_{i=1 \atop (i,p)=1}^{sm-1}\frac{1}{i}\equiv0 \pmod{p^{\alpha}}.
    \end{equation*}
    Therefore, we have
    \begin{equation*}
    \begin{aligned}
        &\sum_{t=0}^{p-1}\sum_{s=1 \atop (s,p)=1}^{np^{\alpha-1}-1}\frac{1}{s^2}\left(1-\frac{2np^{\alpha-1}t}{s}\right)\sum_{i=sm \atop (i,p)=1}^{(s+np^{\alpha-1}t)m-1}\frac{1}{i}\\
        \equiv & \sum_{t=0}^{p-1}\sum_{s=1 \atop (s,p)=1}^{np^{\alpha-1}-1}\frac{1}{s^2}(1-\frac{2np^{\alpha-1}t}{s})\sum_{i=sm \atop (i,p)=1}^{np^{\alpha-1}tm-1}\frac{1}{i}\\
        + &\sum_{t=0}^{p-1}\sum_{s=1 \atop (s,p)=1}^{np^{\alpha-1}-1}\frac{1}{s^2}(1-\frac{2np^{\alpha-1}t}{s})\sum_{i=np^{\alpha-1}tm+1 \atop (i,p)=1}^{(s+np^{\alpha-1}t)m-1}\frac{1}{i} \\
        \equiv &\sum_{t=0}^{p-1}\sum_{s=1 \atop (s,p)=1}^{np^{\alpha-1}-1}\frac{1}{s^2}(1-\frac{2np^{\alpha-1}t}{s})\sum_{i=sm \atop (i,p)=1}^{np^{\alpha-1}tm-1}\frac{1}{i}\\
        + &\sum_{t=0}^{p-1}\sum_{s=1 \atop (s,p)=1}^{np^{\alpha-1}-1}\frac{1}{s^2}(1-\frac{2np^{\alpha-1}t}{s})\sum_{i=1 \atop (i,p)=1}^{sm-1}\frac{1}{i+np^{\alpha-1}tm}\\
        \equiv &\sum_{t=0}^{p-1}\sum_{s=1 \atop (s,p)=1}^{np^{\alpha-1}-1}\frac{1}{s^2}(1-\frac{2np^{\alpha-1}t}{s})\sum_{i=sm \atop (i,p)=1}^{np^{\alpha-1}tm-1}\frac{1}{i}\\
        + &\sum_{t=0}^{p-1}\sum_{s=1 \atop (s,p)=1}^{np^{\alpha-1}-1}\frac{1}{s^2}(1-\frac{2np^{\alpha-1}t}{s})\sum_{i=1 \atop (i,p)=1}^{sm-1}\frac{1}{i}\\
        - &np^{\alpha-1}m\sum_{t=0}^{p-1}t\sum_{s=1 \atop (s,p)=1}^{np^{\alpha-1}-1}\frac{1}{s^2}(1-\frac{2np^{\alpha-1}t}{s})\sum_{i=1 \atop (i,p)=1}^{sm-1}\frac{1}{i^2}\\
        \equiv &\sum_{t=0}^{p-1}\sum_{s=1 \atop (s,p)=1}^{np^{\alpha-1}-1}\frac{1}{s^2}(1-\frac{2np^{\alpha-1}t}{s})\sum_{i=1 \atop (i,p)=1}^{np^{\alpha-1}tm-1}\frac{1}{i}\pmod{p^{\alpha}}.\\
    \end{aligned}
    \end{equation*}
 By Lemma 2.1, we have
 \begin{equation*}
    \sum_{t=0}^{p-1}\sum_{s=1 \atop (s,p)=1}^{np^{\alpha-1}-1}\frac{1}{s^2}(1-\frac{2np^{\alpha-1}t}{s})\sum_{i=sm \atop (i,p)=1}^{(s+np^{\alpha-1}t)m-1}\frac{1}{i}\equiv0\pmod{p^{\alpha}}.
 \end{equation*}
 Hence
 \begin{equation*}
     S(m,np^{\alpha};p)\equiv\sum_{t=0}^{p-1}\sum_{s=1 \atop (s,p)=1}^{np^{\alpha-1}-1}\frac{1}{s^2}\sum_{i=1 \atop (i,p)=1}^{sm-1}\frac{1}{i}\equiv pS(m,np^{\alpha-1};p)\equiv...\equiv p^{\alpha-1}S(m,np;p)\pmod{p^{\alpha}}.
 \end{equation*}
\end{lemma}

\begin{lemma}
    Let p be an odd prime, n, m be positive integers which coprime to p and $k$ be an integer with$k \geq 0$, then we have
    \begin{equation}
      \sum_{a=1 \atop (a,p)=1}^{np-1}\frac{1}{a^k}\left[ \frac{am}{p}\right]\equiv
\begin{cases}
\begin{aligned}
     & -\frac{m+1}{2}n \pmod{p},&&k=0, \\
     & 0\pmod{p},&&1\leq k \leq p-2,\textit{k is even,}\\
     &n\frac{m^k-m^p}{k}B_{p-k}\pmod{p},&&1\leq k \leq p-2,\textit{k is odd},
\end{aligned}
\end{cases}
    \end{equation}
    where $[x]$ denotes the largest integer less than or equal to $x$.\\
    Proof: From Ref. \cite{Cai2016}, we have
    \begin{equation}
        \sum_{a=1}^{p-1}\frac{1}{a^k}\left[ \frac{am}{p}\right]\equiv
\begin{cases}
\begin{aligned}
     & -\frac{m+1}{2} \pmod{p},&&k=0, \\
     & 0\pmod{p},&&1\leq k \leq p-2,\textit{k is even,}\\
     &\frac{m^k-m^p}{k}B_{p-k}\pmod{p},&&1\leq k \leq p-2,\textit{k is odd}.
\end{aligned}
\end{cases}
    \end{equation}
    For any positive integer $b$, we can get
    \begin{equation}
    \begin{aligned}
        \sum_{a=bp+1 \atop (a,p)=1}^{(b+1)p-1}\frac{1}{a^k}\left[ \frac{am}{p}\right]&\equiv\sum_{a=1 }^{p-1}\frac{1}{(a+bp)^k}\left[ \frac{(a+bp)m}{p}\right]\\
        &\equiv\sum_{a=1 \atop }^{p-1}\frac{1}{a^k}\left[ \frac{am}{p}\right]+bm\sum_{a=1}^{p-1}\frac{1}{a^k}\\
        &\equiv\sum_{a=1}^{p-1}\frac{1}{a^k}\left[ \frac{am}{p}\right]\pmod{p}.
    \end{aligned}
    \end{equation}
    So we have proved the lemma.
\end{lemma}
\noindent{\bf Note}.Specially, from Lemma 2.4, we have
\begin{equation}
    T(m,np;p)=\sum_{a=1 \atop (a,p)=1}^{n-1}\frac{1}{a^2}\sum_{i=1 \atop i \equiv sm \pmod{p}}^{am-1}\frac{1}{i}\equiv\frac{1}{m}\sum_{a=1}^{np-1}\frac{1}{a^3}\left[ \frac{am}{p}\right]\equiv\frac{n(m^3-m)}{3m}B_{p-3}\pmod{p}.
\end{equation}

\begin{lemma}
    Let p be an odd prime, n, m be positive integers coprime to p and $\alpha \geq 2$ be a positive integer, then we have
    \begin{equation*}
       T(m,np^{\alpha};p)\equiv \frac{n(m^3-m)}{3m}p^{\alpha-1}B_{p-3}\pmod{p^{\alpha}}.
    \end{equation*}
   Proof: Let $a=s+np^{\alpha-1}t$, $1\leq s \leq np^{\alpha-1}-1$, $0 \leq t \leq p-1$ and $(s,p)=1$, then
   \begin{equation}
       \begin{aligned}
           T(m,np^{\alpha};p)&= \sum_{a=1 \atop (a,p)=1}^{np^{\alpha}-1}\frac{1}{a^2}\sum_{i=1 \atop i\equiv am \pmod{p}}^{am-1}\frac{1}{i}=\sum_{t=0}^{p-1}\sum_{s=1 \atop (s,p)=1}^{np^{\alpha-1}-1}\frac{1}{(s+np^{\alpha-1}t)^2}\sum_{i=1 \atop i \equiv sm \pmod{p}}^{(s+np^{\alpha-1}t)m-1}\frac{1}{i}\\
           &\equiv\sum_{t=0}^{p-1}\sum_{s=1 \atop (s,p)=1}^{np^{\alpha-1}-1}\frac{1}{s^2}(1-\frac{2np^{\alpha-1}t}{s})\left(\sum_{i=1 \atop i \equiv sm \pmod{p}}^{sm-1}\frac{1}{i}+\sum_{i=sm \atop i \equiv sm \pmod{p}}^{(s+np^{\alpha-1}t)m-1}\frac{1}{i}\right)\pmod{p^{\alpha}}.
       \end{aligned}
   \end{equation}
   It is easy to see that
   \begin{equation}
       2np^{\alpha-1}\sum_{t=0}^{p-1}t\sum_{s=1 \atop (s,p)=1}^{p^{\alpha-1}-1}\frac{1}{s^3}\sum_{i=1 \atop i \equiv sm \pmod{p}}^{sm-1}\frac{1}{i}\pmod{p^{\alpha}}.
   \end{equation}
   Since
   \begin{equation}
       \begin{aligned}
       \sum_{i=sm \atop i \equiv sm \pmod{p}}^{(s+np^{\alpha-1}t)m-1}\frac{1}{i} &\equiv\sum_{j=0}^{np^{\alpha-2}tm-1}\frac{1}{sm+jp}\equiv\sum_{j=0}^{np^{\alpha-2}tm-1}\frac{1}{sm\left(1+\frac{jp}{sm}\right)}\equiv  \sum_{j=0}^{np^{\alpha-2}tm-1}\frac{1}{sm}\sum_{k=0}^{\alpha-1}\left(-\frac{jp}{sm}\right)^k\\
       &\equiv\sum_{k=0}^{\alpha-1}\frac{(-p)^k}{(sm)^{k+1}}\sum_{j=0}^{np^{\alpha-2}tm-1}j^k\equiv0\pmod{p^{\alpha-2}}.
       \end{aligned}
   \end{equation}
   With Lemma 2.1, we have
   \begin{equation}
       \sum_{t=0}^{p-1}\sum_{s=1 \atop (s,p)=1}^{np^{\alpha-1}-1}\frac{1}{s^2}(1-\frac{2np^{\alpha-1}t}{s})\sum_{i=sm \atop i \equiv sm \pmod{p}}^{(s+np^{\alpha-1}t)m-1}\frac{1}{i}\equiv0 \pmod{p^{\alpha}}.
   \end{equation}
   Hence by recurrence method, we have
   \begin{equation}
       T(np^{\alpha};p)\equiv pT(np^{\alpha-1};p)\equiv...\equiv p^{\alpha-1}T(np;p)\pmod{p^{\alpha}}.
   \end{equation}
   Finally,we have proved Lemma 2.5.
\end{lemma}

\vskip 8mm

\section{Proof of Theorems }

In this part, we will give the proofs of Theorem 1 and Theorem 2. \\
Proof of Theorem 1. By inclusion-exclusion principle, we have
\begin{equation*}
    \begin{aligned}
        Z(p_1...p_n)=\sum_{i+j+k=p_1...p_n \atop i,j,k \in P_{p_i},1\leq i \leq n}\frac{1}{ijk}&=\sum_{i+j+k=p_1...p_n \atop i,j,k \in P_{p_1}}\frac{1}{ijk}-\sum_{l_1=2}^{n}\sum_{i+j+k=p_1...p_n \atop i,j,k \in P_{p_1},i,j,k\notin P_{p_{l_1}}}\frac{1}{ijk}\\
        &+\sum_{l_1,l_2 \in \{2,...,n\} \atop l_1 \neq l_2}\sum_{i+j+k=p_1...p_n \atop i,j,k \in P_{p_1},i,j,k\notin P_{p_{l_1l_2}}}\frac{1}{ijk}+...\\
        &+(-1)^{n-1}\sum_{i+j+k=p_1...p_n \atop i,j,k \in P_{p_1},i,j,k\notin P_{p_2},P_{p_3},...,P_{p_n}}\frac{1}{ijk}.
    \end{aligned}
\end{equation*}
Then we consider
\begin{equation}
   \begin{aligned}
    \sum_{i+j+k=p_1...p_n \atop i,j,k \in P_{p_1},i,j,k\notin P_{p_{l_1}},P_{p_{l_2}},...,P_{p_{l_n}}}\frac{1}{ijk} &=3\left(\sum_{a=1}^{\frac{p_1...p_n}{p_{l_1}...p_{l_s}}-1}\frac{1}{p_1..p_n-ap_{l_1}...p_{l_s}}\sum_{i+j=ap_{l_1}...p_{l_s} \atop (ij,p_1)=1}\frac{1}{ij} \right. \\
    &\left.-\sum_{a+b+c=\frac{p_1...p_n}{p_{l_1}...p_{l_s}}\atop a,b,c \in P_{p_1}}\frac{1}{ap_{l_1}...p_{l_s}bp_{l_1}...p_{l_s}cp_{l_1}...p_{l_s}}\right)\\
    &+\sum_{a+b+c=\frac{p_1...p_n}{p_{l_1}...p_{l_s}}\atop a,b,c \in P_{p_1}}\frac{1}{ap_{l_1}...p_{l_s}bp_{l_1}...p_{l_s}cp_{l_1}...p_{l_s}}\\
    &=3\sum_{a=1}^{\frac{p_1...p_n}{p_{l_1}...p_{l_s}}-1}\frac{1}{p_1..p_n-ap_{l_1}...p_{l_s}}\sum_{i+j=ap_{l_1}...p_{l_s} \atop (ij,p_1)=1}\frac{1}{ij}\\
    &-2\sum_{a+b+c=\frac{p_1...p_n}{p_{l_1}...p_{l_s}}\atop a,b,c \in P_{p_1}}\frac{1}{ap_{l_1}...p_{l_s}bp_{l_1}...p_{l_s}cp_{l_1}...p_{l_s}}
   \end{aligned}
\end{equation}
And due to Ref.\cite{Cai2016}, we have
\begin{equation*}
    \begin{aligned}
      \sum_{a+b+c=\frac{p_1...p_n}{p_{l_1}...p_{l_s}}\atop a,b,c \in P_{p_1}}\frac{1}{ap_{l_1}...p_{l_s}bp_{l_1}...p_{l_s}cp_{l_1}...p_{l_s}}  \equiv\frac{p_2...p_n}{(p_{l_1}...p_{l_s})^4}Z(p_1)\pmod{p_1}.
    \end{aligned}
\end{equation*}
The first sum in (3.1) equals
\begin{equation*}
    \begin{aligned}
        &\sum_{a=1}^{\frac{p_1...p_n}{p_{l_1}...p_{l_s}}-1}\frac{1}{p_1..p_n-ap_{l_1}...p_{l_s}}\frac{1}{ap_{l_1}...p_{l_s}}\sum_{i+j=ap_{l_1}...p_{l_s} \atop (ij,p_1)=1}\frac{i+j}{ij}\\
        \equiv&-\frac{2}{(p_{l_1}...p_{l_s})^2}\sum_{a=1}^{\frac{p_1...p_n}{p_{l_1}...p_{l_s}}-1}\frac{1}{a^2}\sum_{i=1\atop (i,p_1)=(ap_{l_1}...p_{l_s}-i,p_1)=1}^{ap_{l_1}...p_{l_s}}\frac{1}{i}\\
        \equiv&-\frac{2}{(p_{l_1}...p_{l_s})^2}\sum_{a=1}^{\frac{p_1...p_n}{p_{l_1}...p_{l_s}}-1}\frac{1}{a^2}\sum_{i=1\atop (i,p_1)=1}^{ap_{l_1}...p_{l_s}}\frac{1}{i}+\frac{2}{(p_{l_1}...p_{l_s})^2}\sum_{a=1}^{\frac{p_1...p_n}{p_{l_1}...p_{l_s}}-1}\frac{1}{a^2}\frac{1}{ap_{l_1}...p_{l_s}}\left[\frac{ap_{l_1}...p_{l_s}}{p_1}\right]\\
        \equiv&-\frac{2}{(p_{l_1}...p_{l_s})^2}\sum_{a=1}^{\frac{p_1...p_n}{p_{l_1}...p_{l_s}}-1}\frac{1}{a^2}\sum_{i=1\atop (i,p_1)=1}^{ap_{l_1}...p_{l_s}}\frac{1}{i}+\frac{2}{(p_{l_1}...p_{l_s})^3}\sum_{a=1}^{\frac{p_1...p_n}{p_{l_1}...p_{l_s}}-1}\frac{1}{a^3}\left[\frac{ap_{l_1}...p_{l_s}}{p_1}\right]\\
    \end{aligned}
\end{equation*}
By Lemma 2.2, we have
\begin{equation}
    \sum_{a=1}^{\frac{p_1...p_n}{p_{l_1}...p_{l_s}}-1}\frac{1}{a^2}\sum_{i=1\atop (i,p_1)=1}^{ap_{l_1}...p_{l_s}}\frac{1}{i}\equiv p_2...p_np_{l_1}...p_{l_s}B_{p_1-3}\pmod{p_1}.
\end{equation}
By Lemma 2.4, we have
\begin{equation}
    \sum_{a=1}^{\frac{p_1...p_n}{p_{l_1}...p_{l_s}}-1}\frac{1}{a^3}\left[\frac{ap_{l_1}...p_{l_s}}{p_1}\right]\equiv p_2...p_n\frac{(p_{l_1}...p_{l_s})^3-(p_{l_1}...p_{l_s})}{3}B_{p_1-3}\pmod{p_1}.
\end{equation}
So we have
\begin{equation}
\begin{aligned}
    Z(p_1...p_n)&\equiv-2p_2...p_nB_{p_1-3}\left(1+\sum_{l_1=2}^{n}\frac{2}{p_{l_1}^4}-\sum_{l_1,l_2 \in \{2,...,n\} \atop l_1 \neq l_2}\frac{2}{p_{l_1}^4p_{l_2}^4}+...+(-1)^n\frac{2}{p_{2}^4...p_{n}^4}  \right)\\
    &+6p_2...p_nB_{p_1-3}\left(\sum_{l_1=2}^{n}\frac{2p_{l_1}^2+1}{3p_{l_1}^3}-\sum_{l_1,l_2 \in \{2,...,n\} \atop l_1 \neq l_2}\frac{2p_{l_1}^2p_{l_2}^2+1}{3p_{l_1}^3p_{l_2}^3}+...\right.\\
    &\left.+(-1)^n\frac{2p_{2}^2...p_{n}^2+1}{3p_{2}^3...p_{n}^3} \right)\pmod{p_1}.
\end{aligned}
\end{equation}
We have proved Theorem 1.\\
Proof of Theorem 2. By Lemma 2.3 and Lemma 2.5, we can get
\begin{equation}
\begin{aligned}
    Z(p_1^{r_1}...p_n^{r_n})&\equiv-2p_1^{r_1-1}p_2^{r_2}...p_n^{r_n}B_{p_1-3}\left(1+\sum_{l_1=2}^{n}\frac{2}{p_{l_1}^4}-\sum_{l_1,l_2 \in \{2,...,n\} \atop l_1 \neq l_2}\frac{2}{p_{l_1}^4p_{l_2}^4}+...+(-1)^n\frac{2}{p_{2}^4...p_{n}^4}  \right)\\
    &+6p_1^{r_1-1}p_2^{r_2}...p_n^{r_n}B_{p_1-3}\left(\sum_{l_1=2}^{n}\frac{2p_{l_1}^2+1}{3p_{l_1}^3}-\sum_{l_1,l_2 \in \{2,...,n\} \atop l_1 \neq l_2}\frac{2p_{l_1}^2p_{l_2}^2+1}{3p_{l_1}^3p_{l_2}^3}+...\right.\\
    &\left.+(-1)^n\frac{2p_{2}^2...p_{n}^2+1}{3p_{2}^3...p_{n}^3} \right)\pmod{p_1^{r_1}}.
\end{aligned}
\end{equation}
Therefore, we finally prove Theorem 2.

\hfill
$\square$

\vskip 3mm

\noindent{\bf Acknowledgement}.
This work is supported in part by Shaanxi Fundamental Science Research Project for Mathematics and Physics
(Grant No. 23JSY033) and Natural Science Basic Research Project of Shaanxi Province(2021JM-044).

\end{document}